\documentclass{scrartcl}

\usepackage[english]{babel}
\usepackage[utf8]{inputenc}
\usepackage[T1]{fontenc}
\usepackage{mo}

\makeatletter
\theoremstyle{m@faktstyle}
\newtheorem*{result*}{\IfLanguageName{ngerman}{Satz}{Theorem}}
\makeatother

\usepackage[
    url=false,
    doi=false,
    eprint=false,
    backend=biber
]{biblatex}
\renewcommand*{\biburlsetup}{%
  \Urlmuskip=0mu plus 3mu\relax
  \mathchardef\UrlBigBreakPenalty=100\relax
  \mathchardef\UrlBreakPenalty=200\relax
  \def\UrlBigBreaks{\do\:\do\-}%
  \def\UrlBreaks{%
    \do\.\do\@\do\/\do\\\do\!\do\_\do\|\do\;\do\>\do\]\do\)\do\}%
    \do\,\do\?\do\'\do\+\do\=\do\#\do\$\do\&\do\*\do\^\do\"\do\:}%
  \ifnumgreater{\value{biburlnumpenalty}}{0}
    {\def\do##1{\appto\UrlSpecials{\do##1{\mathchar`##1 \penalty\value{biburlnumpenalty}}}}%
     \do\1\do\2\do\3\do\4\do\5\do\6\do\7\do\8\do\9\do\0}
    {}%
  \ifnumgreater{\value{biburlucpenalty}}{0}
    {\def\do##1{\appto\UrlSpecials{\do##1{\mathchar`##1 \penalty\value{biburlucpenalty}}}}%
     \do\A\do\B\do\C\do\D\do\E\do\F\do\G\do\H\do\I\do\J
     \do\K\do\L\do\M\do\N\do\O\do\P\do\Q\do\R\do\S\do\T
     \do\U\do\V\do\W\do\X\do\Y\do\Z}
    {}%
  \ifnumgreater{\value{biburllcpenalty}}{0}
    {\def\do##1{\appto\UrlSpecials{\do##1{\mathchar`##1 \penalty\value{biburllcpenalty}}}}%
     \do\a\do\b\do\c\do\d\do\e\do\f\do\g\do\h\do\i\do\j
     \do\k\do\l\do\m\do\n\do\o\do\p\do\q\do\r\do\s\do\t
     \do\u\do\v\do\w\do\x\do\y\do\z}
    {}%
  \let\do=\noexpand}

\usepackage{breakurl}

\title{The mod-$2$ cohomology of $32\Gamma_3f$}
\author{Markus Oehme}
\begin{document}
\maketitle

\begin{abstract}
  We compute the group cohomology of $32\Gamma_3f$, a certain group of
  order 32. For this we construct explicit cocycle representatives of
  the cohomology generators. We thus lay to rest a discrepancy between
  several published computations of this cohomology ring.

\end{abstract}

\section{Introduction}

We want to examine the group $32\Gamma_3f$ given by
\begin{equation*}
  \left\<x, y\middle| y^8=1, x^4=y^4, \presup{x}{y}=y^3 \right\> \,.
\end{equation*}
This is a group of size 32 and the notation $32\Gamma_3f$ is taken
from Hall--Senior~\cite{Hall1964}. Alternatively it is the group number
15 of size 32 in the numbering of the Small Groups
library~\cite{Besche2002}. The above presentation is as
in~\cite{Huebschmann1989}.

The mod-$2$-cohomology of $32\Gamma_3f$ has been calculated four times
already, but with two incompatible results. Specifically
Rusin~\cite{Rusin1987} and Huebschmann\cite{Huebschmann1989} state
that all cohomology classes in degree $3$ are nilpotent, whereas
Carlson~\cite{Carlson2003} and Green--King~\cite{Green2011,
  Green2011w} find that there is an indecomposable element in degree 3
which is non-nilpotent. Here we present an independent computation
which verifies the result of Carlson and Green--King; in particular we
exhibit a standard 3-cocycle whose cohomology class is non-nilpotent.

The group admits the following presentation, which we are going to use
throughout the article.
\begin{align*}
  32\Gamma_3f &= \left\<f_1, f_2, f_3, f_4, f_5\middle|\begin{gathered}
      f_1^2 = f_4, f_2^2 = f_3, f_3^2=f_4^2=f_5, f_5^2 = 1, \\
      f_2^{f_1}=f_2f_3, f_3^{f_1}=f_3f_5, f_i^{f_j} = f_i\text{ for all other $j<i$}
    \end{gathered}
  \right\> \\
            &\stackrel{\text{as set}}{=}\left\{f_1^af_2^bf_3^cf_4^df_5^e
              \mid a,b,c,d,e\in\FF_2\right\}
\end{align*}
This compares to the first presentation as follows:
$f_1 \stackrel{\wedge}{=} x$, $f_2 \stackrel{\wedge}{=} y$,
$f_3 \stackrel{\wedge}{=} y^2$, $f_4 \stackrel{\wedge}{=} x^2$ and
$f_5 \stackrel{\wedge}{=} x^4 = y^4$. The operation is then given as
\begin{align*}
  g_1g_2 = (f_1^{a_1}f_2^{b_1}f_3^{c_1}f_4^{d_1}f_5^{e_1})
           (f_1^{a_2}f_2^{b_2}f_3^{c_2}f_4^{d_2}f_5^{e_2})
         = f_1^{a_1+a_2}f_2^{b_1+b_1}f_3^{c_1+c_2+\tilde c}f_4^{d_1+d_2+\tilde d}
           f_5^{e_1+e_2+\tilde e}
\end{align*}
where
\begin{align*}
  \tilde c &= b_1a_2+b_1b_2 \\
  \tilde d &= a_1a_2 \\
  \tilde e &= d_1d_2+a_1a_2d_2+a_1d_1a_2
             +c_1c_2+b_1b_2c_2+b_1c_1b_2+b_1a_2c_2+b_1c_1a_2+c_1a_2+b_1a_2b_2
\end{align*}
are the deviations from the elementary abelian case.

Along with~\cite[Conjecture 3]{Eick2015} we propose the following
conjecture.

\begin{verm*}
  The family of groups
  $\Phi_n := \left\<x, y\middle| y^{2^n}=1, x^4=y^{2^{n-1}},
    \presup{x}{y}=y^{2^n-1} \right\> \cong \left\<x, y\middle|
    y^{2^n}=1, x^4=y^{2^{n-1}}, \presup{x}{y}=y^{2^{n-1}-1} \right\>$
  for $n\ge3$ is the only exception to Huebschmanns results and all
  groups of the family have the cohomology
  \[
  H^*(\Phi_n) \cong \FF_2[u_1, v_1, W_2, y_3, z_4]/(v^2+uv = u^2 = uW
  = y^2 + W^3 = uy = 0)\,.
  \]
\end{verm*}
Thus case (2) of Theorem~E of~\cite{Huebschmann1989} splits into two
cases, the $\Phi_n$ with the corrected cohomology and the remainder
with the already correct cohomology.

\section{Breaking the problem apart}

We describe $32\Gamma_3f$ as a central extension
\begin{equation}
  1\to C_2=\<f_5\> \hookrightarrow 32\Gamma_3f \twoheadrightarrow 16\Gamma_2c_2 \to 1
  \label{eq:extension32}\tag{I}
\end{equation}
where the quotient is the group $16\Gamma_2c_2$ of size 16 according
to Hall--Senior or the group number 4 of size 16 according to the Small
Groups Library. We describe $16\Gamma_2c_2$ as
\[
  16\Gamma_2c_2 = \left\<f_1, f_2, f_3, f_4, f_5\middle| f_5 = 1\right\>
\]
where the $f_i$ otherwise behave as in the case of
$32\Gamma_3f$. Furthermore we describe $16\Gamma_2c_2$ as a central
extension
\begin{equation}
  1\to C_2=\<f_4\> \hookrightarrow 16\Gamma_2c_2 \twoheadrightarrow D_8 \to 1
  \label{eq:extension16}\tag{II}
\end{equation}
where
\[
  D_8 = \left\<f_1, f_2, f_3, f_4, f_5\middle| f_4 = f_5 = 1\right\>
\]
is easily identifiable as a dihedral group of order 8.

\section{Cohomology}

We denote by $\bar a$, $\bar b$, $\bar c$, $\bar d$, $\bar e$, the set
maps $32\Gamma_3f\to\FF_2$ mapping $f_1^af_2^bf_3^cf_4^df_5^e$ to the
respective exponent. We now utilize the bar construction to describe
the cohomology of the group $G$ ($32\Gamma_3f$ or a quotient thereof),
so that in degree $n$ the cochains are given by the vector space of
maps $G^n\to\FF_2$ as given above, but with indices to indicate the
relevant copy of $G$. Thus for example
$\bar b_1\bar a_2\bar c_2:G^2\to\FF_2:(g_1, g_2) =
(f_1^{a_1}f_2^{b_1}f_3^{c_1}f_4^{d_1}f_5^{e_1},
f_1^{a_2}f_2^{b_2}f_3^{c_2}f_4^{d_2}f_5^{e_2})\mapsto b_1a_2c_2$.

The cohomology of the $D_8$ is (see for example~\cite[Theorem
2.7]{Adem2004}, where the group generators are different so that our
$u$ is the sum of their degree one cohomology generators)
\[
H^*(D_8) = \FF_2[u_1, v_1, w_2]/(v^2+uv = 0)
\]
where the classes of degree one are given by $u = [\bar a]$,
$v = [\bar b]$ with square brackets denoting the appropriate
equivalence classes and the subscripts denoting the degree. Further we
calculate for $w$ in degree two (as $\<f_2, f_3\>$ is a $C_4$ we have
$w\approx [\bar c^2]$).
\begin{align*}
  \delta({\bar c}_1{\bar c}_2&+{\bar b}_1{\bar a}_2{\bar c}_2+{\bar b}_1{\bar c}_1{\bar a}_2+{\bar b}_1{\bar b}_2{\bar c}_2+{\bar b}_1{\bar c}_1{\bar b}_2+{\bar c}_1{\bar a}_2+{\bar b}_1{\bar a}_2{\bar b}_2)(g_1, g_2, g_3)  =: \delta(\Xi)(g_1, g_2, g_3)\\
               &=\Xi(g_2, g_3) - \Xi(g_1g_2, g_3) + \Xi(g_1, g_2g_3) - \Xi(g_1, g_2) \\
               &=(c_2c_3+b_2a_3c_3+b_2c_2a_3+b_2b_3c_3+b_2c_2b_3+c_2a_3+b_2a_3b_3) \\
               &\qquad - \bigl((c_1+c_2+b_1a_2+b_1b_2)c_3+(b_1+b_2)a_3c_3+(b_1+b_2)(c_1+c_2+b_1a_2+b_1b_2)a_3 \\
               &\qquad\qquad +(b_1+b_2)b_3c_3+(b_1+b_2)(c_1+c_2+b_1a_2+b_1b_2)b_3\\
               &\qquad\qquad +(c_1+c_2+b_1a_2+b_1b_2)a_3+(b_1+b_2)a_3b_3\bigr) \\
               &\qquad + \bigl(c_1(c_2+c_3+b_2a_3+b_2b_3)+b_1(a_2+a_3)(c_2+c_3+b_2a_3+b_2b_3)+b_1c_1(a_2+a_3) \\
               &\qquad\qquad +b_1(b_2+b_3)(c_2+c_3+b_2a_3+b_2b_3)+b_1c_1(b_2+b_3)+c_1(a_2+a_3) \\
               &\qquad\qquad +b_1(a_2+a_3)(b_2+b_3)\bigr) \\
               &\qquad - (c_1c_2+b_1a_2c_2+b_1c_1a_2+b_1b_2c_2+b_1c_1b_2+c_1a_2+b_1a_2b_2) \\
               &= 0.
\end{align*}
Hence $w = [{\bar c}_1{\bar c}_2+{\bar b}_1{\bar a}_2{\bar c}_2+{\bar b}_1{\bar c}_1{\bar a}_2+{\bar b}_1{\bar b}_2{\bar c}_2+{\bar b}_1{\bar c}_1{\bar b}_2+{\bar c}_1{\bar a}_2+{\bar b}_1{\bar a}_2{\bar b}_2]$.

Furthermore the only non-trivial Steenrod operation is $Sq^1(w) = uw$
(see again~\cite[Theorem 2.7]{Adem2004}).

Now we can determine the cocycle of the
extension~(\ref{eq:extension16}). To do this we select a splitting
$\sigma: D_8\to 16\Gamma_2c_2$ on the set level and compute
$q:D_8^2\to \<f_4\>:(g_1, g_2)\mapsto\sigma(g_1) \sigma(g_2)
\sigma(g_1g_2)^{-1}$,
which is a representative for the cocycle in $H^2(D_8, \<f_4\>\cong\FF_2)$. We
choose the canonical splitting mapping $f_1^af_2^bf_3^c$ to
$f_1^af_2^bf_3^cf_4^0$. Thus we get
\begin{align*}
q(g_1, g_2)& = q(f_1^{a_1}f_2^{b_1}f_3^{c_1}, f_1^{a_2}f_2^{b_2}f_3^{c_2}) \\
           & = \sigma(f_1^{a_1}f_2^{b_1}f_3^{c_1}) \sigma(f_1^{a_2}f_2^{b_2}f_3^{c_2})
             \sigma(f_1^{a_1+a_2}f_2^{b_1+b_1}f_3^{c_1+c_2+b_1a_2+b_1b_2})^{-1} \\
           & = f_1^{a_1}f_2^{b_1}f_3^{c_1}f_1^{a_2}f_2^{b_2}f_3^{c_2}
             (f_1^{a_1+a_2}f_2^{b_1+b_1}f_3^{c_1+c_2+b_1a_2+b_1b_2})^{-1} \\
           & = f_1^{a_1+a_2}f_2^{b_1+b_1}f_3^{c_1+c_2+b_1a_2+b_1b_2}f_4^{a_1a_2}
             (f_1^{a_1+a_2}f_2^{b_1+b_1}f_3^{c_1+c_2+b_1a_2+b_1b_2})^{-1} \\
           & = f_1^{a_1+a_2}f_2^{b_1+b_1}f_3^{c_1+c_2+b_1a_2+b_1b_2}
             (f_1^{a_1+a_2}f_2^{b_1+b_1}f_3^{c_1+c_2+b_1a_2+b_1b_2})^{-1}f_4^{a_1a_2} \\
           & = f_4^{a_1a_2}
\end{align*}
and the cocycle $\alpha$ is $[\bar a]^2 = u^2$.

\begin{table}[h!]
  \centering
  \centerline{\begin{tabular}{c||c|c|c|c|c|c}
  5 & $\vdots{}$  &&&&&\\ \hline
  2 & $t^2$  & $\hdots{}$ &&&&\\ \hline
  1 & $t$  & $ut$, $vt$ & $\hdots{}$ &&&\\ \hline
  0 & 1 & $u$, $v$ & $u^2$, $uv=v^2$, $w$ & $u^3$, $v^3$, $uw$, $vw$ & $u^4$, $v^4$, $u^2w$, $w^2$, $v^2w$ & $\cdots{}$\\ \hline\hline
  ~ & 0 & 1 & 2 & 3 & 4
  \end{tabular}
}
\end{table}

Now we use the Leray-Serre-spectral-sequence to determine the
cohomology. We have $E_2^{pq} = H^p(D_8)\otimes H^q(\<f_4\>) \implies
H^*(16\Gamma_2c_2)$, where $H^*(\<f_4\>) = H^*(C_2) = \FF_2[t_1]$ (and
$t = [\bar d]$). Now the cocycle $u^2\in H^2(D_8)$ has to be killed and
hence $d_2(t) = u^2$. Now by the algebra structure we can compute
$d_2$ everywhere and get $E_3 \cong (H^*(D_8)/(u^2))\otimes
\FF_2[t^2]$. Now we use the Steenrod operations to see by Kudo's
transgression theorem
\begin{equation*}
  d_3(t^2) = d_3(Sq^1(t)) = Sq^1(d_2(t)) = Sq^1(u^2) = 0\,.
\end{equation*}
By the algebra structure all of $d_3$ is zero. Furthermore $d_n(t^2)$
for $n>3$ vanishes trivially. Thus we have a collapse at the $E_3$
page and retrieve
\[
H^*(16\Gamma_2c_2) \cong \FF_2[u_1, v_1, w_2, x_2]/(v^2+uv = u^2 = 0)\,.
\]
Here $u$, $v$ and $w$ are simply inflated from $H^*(D_8)$.

Now $x$ is roughly $[\bar d^2]$ and we verify
\begin{align*}
  \delta({\bar d}_1{\bar d}_2&+{\bar a}_1{\bar a}_2{\bar d}_2+{\bar a}_1{\bar d}_1{\bar a}_2)(g_1, g_2, g_3) \\
               &=(d_2d_3+a_2a_3d_3+a_2d_2a_3) \\
               &\qquad - ((d_1+d_2+a_1a_2)d_3+(a_1+a_2)a_3d_3+(a_1+a_2)(d_1+d_2+a_1a_2)a_3) \\
               &\qquad + (d_1(d_2+d_3+a_2a_3)+a_1(a_2+a_3)(d_2+d_3+a_2a_3)+a_1d_1(a_2+a_3)) \\
               &\qquad - (d_1d_2+a_1a_2d_2+a_1d_1a_2) \\
               &= 0.
\end{align*}
Thus
$x=[{\bar d}_1{\bar d}_2+{\bar a}_1{\bar a}_2{\bar d}_2+{\bar
  a}_1{\bar d}_1{\bar a}_2]$.
Furthermore the non-trivial Steenrod operations are $Sq^1(w)=uw$
(still valid) and $Sq^1(x)=0$. For the latter we have to do a bit of
work, the result lies in degree three and hence is a linear
combination of $uw$, $vw$, $ux$ and $vx$. Now we restrict to the two
possible $C_4\times C_2$ subgroups. For $\<f_1, f_3, f_4\>$ with cohomology
$\FF_2[p_1, q_1, r_2]/(q^2=0)$ we get the restrictions
\begin{gather*}
  u \mapsto q, \qquad v \mapsto 0, \qquad w \mapsto p^2+pq, \qquad x \mapsto r
\end{gather*} and for $\<f_2, f_3, f_4\>$ with cohomology
$\FF_2[p_1', q_1', r_2']/(q'^2=0)$ we get the restrictions
\begin{gather*}
  u \mapsto 0, \qquad v \mapsto q', \qquad w \mapsto r', \qquad x \mapsto p'^2.
\end{gather*}
Now $Sq^1(p'^2)=0=Sq^1(r)$ (for the second one see
e.\,g.~\cite[\S\,7.4]{Carlson2003}). Hence $Sq^1(x)$ cannot contain
any of the four aforementioned constituents and thus must vanish.

\subsection{Final computation for $32\Gamma_3f$}

We determine the extension cocycle of~(\ref{eq:extension32}) as
previously.
\begin{align*}
q(g_1, g_2)& = \sigma(f_1^{a_1}f_2^{b_1}f_3^{c_1}f_4^{d_1})
             \sigma(f_1^{a_2}f_2^{b_2}f_3^{c_2}f_4^{d_2})
             \sigma(f_1^{a_1+a_2}f_2^{b_1+b_1}f_3^{c_1+c_2+b_1a_2+b_1b_2}f_4^{d_1+d_2+a_1a_2})^{-1} \\
           & = f_5^{d_1d_2+a_1a_2d_2+a_1d_1a_2+c_1c_2+b_1b_2c_2+b_1c_1b_2+b_1a_2c_2+b_1c_1a_2+c_1a_2+b_1a_2b_2}
\end{align*}
Thus we receive for the cocycle $\alpha$ now the value $w+x$.

\begin{table}[h!]
  \centering
  \centerline{\begin{tabular}{c||c|c|c|c|c|c}
  5 & $\vdots{}$  &&&&&\\ \hline
  4 & $t^4$  & $\hdots{}$ &&&&\\ \hline
  3 & $t^3$  & $\hdots{}$ &&&&\\ \hline
  2 & $t^2$  & $ut^2$, $vt^2$ & $uvt^2=v^2t^2$, $wt^2$, $xt^2$ & $\hdots{}$ &&\\ \hline
  1 & $t$  & $ut$, $vt$ & $uvt=v^2t$, $wt$, $xt$ & $\hdots{}$ &&\\ \hline
  0 & 1 & $u$, $v$ & $uv=v^2$, $w$, $x$ & $uw$, $vw$, $ux$, $wx$ & $v^2w$, $v^2x$, $w^2$, $wx$, $x^2$ & $\cdots{}$\\ \hline\hline
  ~ & 0 & 1 & 2 & 3 & 4 & 5
  \end{tabular}
  }
\end{table}

Again we use the Leray-Serre-spectral-sequence to determine the
cohomology. We have $E_2^{pq} = H^p(16\Gamma_2c_2)\otimes
H^q(\<f_5\>) \implies H^*(32\Gamma_3f)$, where $H^*(\<f_5\>) = H^*(C_2)
= \FF_2[t_1]$ (and $t = [\bar e]$). Now the cocycle $w+x\in
H^2(16\Gamma_2c_2)$ has to be killed and hence $d_2(t) = w+x$. Now by
the algebra structure we can compute $d_2$ everywhere and get $E_3
\cong (H^*(16\Gamma_2c_2)/(w+x))\otimes \FF_2[t^2]$. Again we use
Steenrod operations to determine $d_3$ as
\begin{equation*}
  d_3(t^2) = d_3(Sq^1(t)) = Sq^1(d_2(t)) = Sq^1(w+x) = uw\,.
\end{equation*}
We can thus determine $d_3$ everywhere. Note that $d_3(ut^2) = u^2w =
0$ and $d_3(v^2t^2) = v^2uw = vu^2w = 0$. This gives us $E_3 \cong
(H^*(16\Gamma_2c_2)/(w+x, uw))\otimes \FF_2[t^4]$. Again we use
Steenrod operations for the next possible non-trivial differential
\begin{equation*}
  d_5(t^4) = d_5(Sq^2(t^2)) = Sq^2(d_3(t^2)) = Sq^2(uw) = u^3w + uw^2\, \stackrel{\wedge}{=} 0.
\end{equation*}
As before we see, that the spectral sequence collapses. We end up with
the following $E_\infty$-page, where we omitted equal things due to
$w = x$ and $uv = v^2$.

\begin{table}[h!]
  \centering
  \centerline{\begin{tabular}{c||c|c|c|c|c|c|c|c}
  5 & $\vdots{}$  & $\vdots{}$ & $\vdots{}$ &&&&&\\ \hline
  4 & $t^4$  & $ut^4$, $vt^4$ &  $wt^4$, $v^2t^4$, $u^2t^4=0$
             & $\hdots{}$ &&&&\\ \hline
  3 &&&&&&&&\\ \hline
  2 & & $ut^2$ & $v^2t^2$, $u^2t^2=0$ & $wut^2$ & $vwut^2$ & $\hdots{}$ &&\\ \hline
  1 &&&&&&&&\\ \hline
  0 & 1 & $u$, $v$ & $v^2$, $w$ & $vw$ & $w^2$ & $vw^2$ & $w^3$ & $\cdots{}$\\ \hline\hline
  ~ & 0 & 1 & 2 & 3 & 4 & 5 & 6 & 7
  \end{tabular}
  }
\end{table}

Now we have the generators $u_1$, $v_1$, $W_2$ ($w$ equalling $x$),
$y_3$ given by $ut^2$ and $z_4$ given by $t^4$. We also get the
relations for $u$, $v$ and $W$ and furthermore $z$ has no
relations. However the relations involving $y$ still need to be
determined, since we have $y^2\stackrel{\wedge}{=}u^2t^4 = 0$
and $uy\stackrel{\wedge}{=}u^2t^2 = 0$.

Now $y$ corresponds roughly to $\bar a\bar e^2$. Indeed with
$x_{12} := {\bar d}_1{\bar d}_2+{\bar a}_1{\bar a}_2{\bar d}_2+{\bar
  a}_1{\bar d}_1{\bar a}_2$
and
$w_{12} := {\bar c}_1{\bar c}_2+{\bar b}_1{\bar a}_2{\bar c}_2+{\bar
  b}_1{\bar c}_1{\bar a}_2+{\bar b}_1{\bar b}_2{\bar c}_2+{\bar
  b}_1{\bar c}_1{\bar b}_2+{\bar c}_1{\bar a}_2+{\bar b}_1{\bar
  a}_2{\bar b}_2$ one can verify that the differential of
\medskip
\begin{gather*}
  \Bigl(x_{12}w_{12} + (\bar e_1 + \bar e_2)(x_{12} + w_{12}) + \bar
  b_1\bar c_1\bar a_2 + \bar b_1\bar c_1\bar a_2\bar c_2 + \bar
  b_1\bar d_2 + \bar c_1\bar a_2\bar c_2 + \bar e_1\bar a_2 + \bar
  e_1\bar e_2\Bigr)\bar a_3 \\
  + x_{12}\bar d_3
\end{gather*}
vanishes. Thus we can use this term as representative for $y$.

To determine the relations of $y$ we use the subgroup
$K := \< f_2, f_3, f_4, f_5 \> \cong C_8\times C_2$, where the
generators are $f_2$ and $f_3f_4$. It has cohomology
$H^*(K) \cong \FF_2[\xi_1, \phi_1, \chi_2]/(\phi^2)$, where $\xi$
belongs to the $C_2$.

We first compute the restrictions on the chain level. The inclusion is
$K\hookrightarrow G : f_2^if_3^jf_5^k(f_3f_4)^l \mapsto
f_2^if_3^{j+l}f_4^lf_5^k$ giving restrictions
\begin{gather*}
 \bar a \mapsto 0, \qquad
 \bar b \mapsto \bar i, \qquad
 \bar c \mapsto \bar j + \bar l, \qquad
 \bar d \mapsto \bar l, \qquad
 \bar e \mapsto \bar k.
\end{gather*}
The cohomology classes of further interest are $\xi = [\bar l]$ and
$\phi = [\bar i]$. With them we get restrictions (suppressing most
terms which map to zero because they contain $\bar a$)
\begin{align*}
 u &= [\bar a] \mapsto 0, \\
 v &= [\bar b] \mapsto [\bar i] = \phi, \\
 W &= [\bar d^2 + \bar a(\ldots)] \mapsto [\bar l^2] = \xi^2, \\
 y &= [\bar d^3 + \bar a(\ldots)] \mapsto [\bar l^3] = \xi^3. \\
\end{align*}

Now we look at the spectral sequence and see that during ungrading we
get the following uncertainties: $y^2$ is a linear combination of
$W^3$ and $vWy$ whereas $uy$ is a scalar multiple of $W^2$.

We restrict to $K$ and immediately retrieve $y^2 = W^3$ and $uy =
0$. Thus we get the following.

\begin{result*}
  The mod-$2$-cohomology ring of $32\Gamma_3f$ is
  \[
  H^*(32\Gamma_3f, \FF_2) \cong \FF_2[u_1, v_1, W_2, y_3, z_4]/(v^2+uv = u^2
  = uW = y^2 + W^3 = uy = 0)\,.
  \]
\end{result*}

\printbibliography

\end{document}